\documentclass[12pt]{article}
\usepackage{latexsym,amsfonts,amssymb,mathrsfs,float}
\usepackage{dsfont}
\usepackage{amsthm}
\usepackage[numbers]{natbib}
\usepackage{amsmath}
\usepackage{indentfirst}
\usepackage{bbm}

\usepackage{color}
\allowdisplaybreaks

\newtheorem{thm}{Theorem}[section]

\newtheorem{lem}[thm]{Lemma}
\newtheorem{pro}[thm]{Proposition}

\newtheorem{rmk}[thm]{Remark}

\numberwithin{equation}{section}

\begin{document}

\title{\bf  Limit behavior of linearly edge-reinforced random walks on the half-line}

\author{Zechun Hu$^{1}$,  Renming Song$^{2}$, and Li Wang$^1$\\ \\
  {\small $^1$ College of Mathematics, Sichuan  University,
 Chengdu 610065, China}\\
 {\small zchu@scu.edu.cn; wangli0@stu.scu.edu.cn}\\ \\
  {\small $^2$Department of Mathematics,
University of Illinois Urbana-Champaign, Urbana, IL 61801, USA}\\
 {\small rsong@illinois.edu}}
\date{}
\maketitle

\begin{abstract}

Motivated by the article [M. Takei, Electron. J. Probab. 26 (2021), article no. 104], we study the limit behavior of linearly edge-reinforced random walks on the half-line
$\mathbb{Z}_+$ with reinforcement parameter $\delta>0$, and each edge $\{x,x+1\}$ has the initial weight $x^{\alpha}\ln^{\beta}x$ for $x > 1$ and $1$ for $x = 0, 1$.
The aim of this paper is to study the almost sure limit behavior of the walk in the recurrent regime, and extend the results of Takei mentioned above.
 \end{abstract}


\noindent {\bf 2020 MR Subject Classification (2020)}\quad 60K35

\section{ Introduction}

Linearly edge-reinforced random walks (LERRW),  introduced by Coppersmith and Diaconis \cite{CD}, have attracted the attention of a lot researchers in the last 40 years.
An LERRW  on a connected  graph can be defined as follows.
Every edge is assigned a positive initial  weight. The random walker starts at a vertex
$\mathbf{0}$. In each step, the walker moves to a nearest neighbor by traversing an edge,  with a probability proportional to the weight of that edge. Each time an edge is traversed, its weight is increased by a fixed constant $\delta>0$
({\it linear edge-reinforcement parameter}). The parameter $\delta$ is taken to be 1 in the early studies.

On all locally finite graphs, an LERRW is equivalent to a random walk in a random environment (RWRE), or more precisely, a mixture of Markov chains \cite{DF, MR2}.
This equivalence enhances the tractability of the LERRW and has, for example, led to applications in Bayesian statistics \cite{BS, DR, MR1}.
For acyclic graphs, the recurrence/transience phase transition for this process was first observed by Pemantle on the infinite binary tree \cite{PR}.
Akahori et al. \cite{AJ} obtained some results on the phase transition of LERRW on the half-line, characterizing the trade-off between  the initial weights and reinforcement parameter.  On critical Galton-Watson trees, Andriopoulos and Archer \cite{AG} established an invariance principle in the recurrent regime and an upper bound on the displacement
in the transient regime.

Another interesting topic is on the limit behavior of the LERRW. For initially fair, sequence-type reinforced random walks on the line, Davis \cite{DB} proved the following strong law of large numbers:
$$
\lim_{n\to\infty}\frac{X_n}{n}=0\ \text{a.s.}.
$$
See Takeshima \cite{TM} for a generalization.
The earliest limit theorem  for LERRWs on $b$-ary trees was obtained by Collevecchio \cite{CA}.
He proved that, for $b \geq 70$, a strong law of large numbers holds for the distance of the process from the root, and also established a central limit theorem for sufficiently large $b$. Further related limit theorems for random walks in random
environments
on trees can be found in Aid\'ekon \cite{Aid} and Collevecchio et al. \cite{CTU}.
 Later, Takei \cite{TMA} investigated the almost sure behavior of LERRW in the recurrent regime and obtained a limit theorem, which is a counterpart of the law of the iterated logarithm for simple random walks (see Theorems \ref{thm-2.2} and \ref{thm-2.3} below).

 Compared to LERRWs on acyclic graphs, the study of LERRWs on cyclic graphs is much more challenging.
A key problem, the recurrence/transience of the LERRW on $\mathbb{Z}^{d}$, was posed by Diaconis in the late 1980s.
Except for the one-dimensional case ($d=1$), this problem remained open for many years until breakthroughs were achieved in the past decade.
Key contributions in this area include the following. Angel et al. \cite{AO}, and  Sabot and
Tarr\`{e}s \cite{SC},
 independently and using different approaches, proved that the LERRW on $\mathbb{Z}^{d}$ is recurrent for sufficiently large reinforcement.
Disertori et al. \cite{DM} showed that for dimensions $d \geq 3$ and small reinforcement, LERRW is transient. Sabot and Zeng \cite{SC1} completely resolved the two-dimensional case, proving that the LERRW on $\mathbb{Z}^2$ is recurrent for all initial constant weights; another significant contribution of their work was the establishment of a 0-1 law for transience/recurrence on $\mathbb{Z}^{d}$.
 Furthermore, for LERRWs on $\mathbb{Z}^d$, $d\geq 3$, Poudevigne-Auboiron \cite{PA} proved a monotonicity result with respect to the initial edge weights: increasing the initial weights makes the walk more transient.

In \cite{TMA}, Takei studied the limit behavior of LERRWs on $\mathbb{Z}_+$ in the recurrent regime, where the initial weight of the edge $\{x,x+1\}$, $x\in \mathbb{Z}_+$, is $x^{\alpha}$ for $x \in \mathbb{N}$ and 1 for $x=0$. Motivated by \cite{TMA}, we study in this paper the limit behavior of
LERRWs on $\mathbb{Z}_+$ in the recurrent regime, where the initial weight of the edge $\{x,x+1\}$, $x\in \mathbb{Z}_+$, is $x^{\alpha}\ln^{\beta}x$ for $x=2, 3, \dots$ and 1 for $x=0, 1$.

Our results show that  the term $\ln^{\beta}x$ influences the limit behavior of
the LERRW in the case $\alpha<1$ for any $\beta\neq 0$, and in the case $\alpha=1$ for $\beta<0$.
The arguments of this paper follow the general strategy of \cite{TMA}, with some modifications.

\section{Main results}

We now formally define a linearly edge-reinforced random walk (LERRW) $\boldsymbol{X}=\{X_n\}_{n\ge 0}$ on $\mathbb{Z}_+$.
Let $(w_0(x):  x\in \mathbb{Z}_+)$ be a collection of positive numbers.  $w_0(x)$ stands for the initial weight of the edge $\{x, x+1\}$. Let $\delta>0$ be the reinforcement parameter. For any $n\in \mathbb{N}$
and $x\in\mathbb{Z}_+$, we define $\phi_n(x):=\sum_{i=1}^{n} \mathbf{1}_{\left\{\left\{X_{i-1}, X_{i}\right\}=\{x,x+1\}\right\}}$ to be the number of times the edge $\{x, x+1\}$ has been traversed up to time $n$.  For any $x\in \mathbb{Z}_+$ and $n\in \mathbb{N}$, we define
\begin{align}
w_{n}(x)=w_0(x)+\delta\phi_n(x).
\end{align}
Set $w_{n}(-1)=0$  for all $n\in \mathbb{Z}_+$.
An LERRW
on $\mathbb{Z}_+$ with initial weights
$(w_0(x): x\in \mathbb{Z}_+)$ and reinforcement parameter $\delta$ is a Markov chain on $\mathbb{Z}_+$ with transition probability
\begin{align}
P\left(X_{n+1}=X_{n}+1 \mid X_{0}, \ldots, X_{n}\right) & =1-P\left(X_{n+1}=X_{n}-1 \mid X_{0}, \ldots, X_{n}\right) \nonumber \\
& =\frac{w_{n}\left(X_{n}\right)}{w_{n}\left(X_{n}-1\right)+w_{n}\left(X_{n}\right)}.
\end{align}
The assumption  $w_{n}(-1)=0$  for all $n\in \mathbb{Z}_+$ implies that the origin is a reflection point.
When $\delta=0$, we get a (un-reinforced) random walk.

An LERRW $\boldsymbol{X}$ is said to be recurrent if every point is visited infinitely often, and
 transient otherwise.
Takeshima \cite{TM} proved a recurrence/transience dichotomy for
the LERRW on $\mathbb{Z}_{+}$.
(Takeshima worked with the case $\delta=1$, but his argument works for any $\delta\ge 0$.)  Takei \cite{TMA}  gave a short elementary proof of this dichotomy which is Theorem \ref{thm-2.1} below.

\begin{thm}\label{thm-2.1} {\rm (\cite[Theorem 4.1]{TM})}
 Let $ \boldsymbol{X}$ be an LERRW
 on $ \mathbb{Z}_{+}$ with initial weights
$(w_0(x): x\in \mathbb{Z}_+)$ and reinforcement parameter $\delta$. Let
$\Phi_{0}:=\sum_{x=0}^{\infty}\frac{1}{w_0(x)}.$\\
\indent (i) If  $\Phi_{0}=+\infty$, then $\boldsymbol{X}$ is recurrent a.s..\\
\indent (ii) If  $\Phi_{0}<+\infty$, then $\boldsymbol{X}$ is transient a.s..
\end{thm}

\subsection{Takei's results}

Takei \cite{TMA} studied the LERRW $\boldsymbol{X}$ on $\mathbb{Z}_{+}$ with the initial weights
\begin{align}\label{2.0}
w_0(x)
=\begin{cases}
1 & (x=0),\\
x^{\alpha} & (x \in \mathbb{N}=\{1,2,3,\cdots \}),
\end{cases}
\end{align}
and reinforcement parameter $\delta>0$.
By Theorem \ref{thm-2.1}, $\boldsymbol{X}$ is recurrent a.s. if and only if $\alpha \leq 1$. Ikenami \cite{IM} proved that,  if  $0\leq\alpha<1$ and $\delta=1$, then for any  $\varepsilon>0$,
$\lim\limits_{n \rightarrow \infty} \frac{X_{n}}{(\ln n)^{(1+\varepsilon) /(1-\alpha)}}=0$ a.s..\\
\indent  Takei \cite{TMA} showed that,  if $\alpha<1$, the precise order of oscillation of $X_{n}$  is $(\ln n)^{1 /(1-\alpha)}$.  Let
\begin{align}\label{2.4}
K(\alpha, \delta):=\begin{cases}
\frac{1-\alpha}{2 \delta} & (\alpha<0), \\
\left(\Psi \left(\frac{1}{2 \delta}+\frac{1}{2}\right)-\Psi \left(\frac{1}{2\delta}\right)\right)^{-1} & (\alpha=0), \\
\frac{1-\alpha}{\delta} & (0<\alpha<1),
\end{cases}
\end{align}
where $\Psi(z)=\Gamma^{\prime}(z) / \Gamma(z)$ is the digamma function.

\begin{thm}\label{thm-2.2} {\rm (\cite[Theorem 2.2]{TMA})} Assume that $\alpha<1$ and $\delta>0$.
The LERRW $\boldsymbol{X}$ on $\mathbb{Z}_+$ with initial weights \eqref{2.0} and  reinforcement parameter $\delta$ satisfies that
\begin{align*}
\limsup _{n \rightarrow \infty} \frac{X_{n}}{\{K(\alpha, \delta) \ln n\}^{1 /(1-\alpha)}}=1 \quad \text { a.s., }
\end{align*}
where $K(\alpha, \delta)$ is defined as in \eqref{2.4}.
\end{thm}

\begin{thm}\label{thm-2.3} {\rm (\cite[ Theorem 2.4]{TMA})}
Assume that $\alpha = 1$ and $\delta > 0$, and consider the LERRW $\boldsymbol{X}$ on $\mathbb{Z}_+$ with initial weights \eqref{2.0} and reinforcement parameter $\delta$.\\
\indent (i) If $\delta > 2$, then for any $\varepsilon > 0$,
$$ \limsup_{n \to \infty} \frac{X_n}{n^{(1-\varepsilon)/\delta}} = +\infty
\quad \text{and} \quad  \lim_{n \to \infty} \frac{X_n}{n^{(1+\varepsilon)/\delta}} = 0 \quad  \text {a.s.}.
$$
\indent (ii) If $0 < \delta \leq 2$, then for any $\varepsilon > 0$,
$$\limsup_{n \to \infty} \frac{X_n}{n^{(1-\varepsilon)/2}} = +\infty
\quad \text{and} \quad  \lim_{n \to \infty} \frac{X_n}{n^{(1+\varepsilon)/2}} = 0 \quad  \text {a.s.. }$$
\end{thm}

To see
the effect of the reinforcement, Takei also proved the following
result.

\begin{thm}\label{thm-2.4} {\rm (\cite[Theorem 2.5]{TMA})} Consider
the un-reinforced random walk $\boldsymbol{X}$ on $\mathbb{Z}_+$ with  initial weights \eqref{2.0}.\\
\indent (i) If $\alpha < -1$, then for any $\varepsilon > 0$,
$$\limsup_{n \to \infty} \frac{X_n}{n^{1/(1-\alpha)}} >0
\quad \text{and} \quad \lim_{n \to \infty} \frac{X_n}{\{n(\ln n)^{1+\varepsilon}\}^{1/(1-\alpha)}} = 0 \quad \text{a.s..}
$$
\indent (ii) If $\alpha=-1$, then for any $\varepsilon > 0$,
$$\limsup_{n \to \infty} \frac{X_n}{n^{(1-\varepsilon)/2}} = +\infty
\quad \text{and} \quad \lim_{n \to \infty} \frac{X_n}{\{n(\ln n)^{1+\varepsilon}\}^{1/2}} = 0 \quad \text{a.s..}$$
\indent (iii) If $-1 < \alpha \leq 1$, then for any $\varepsilon > 0$,
$$\limsup_{n \to \infty} \frac{X_n}{n^{1/2}} > 0
\quad \text{and} \quad \lim_{n \to \infty} \frac{X_n}{\{n(\ln n)^{1+\varepsilon}\}^{1/2}} = 0 \quad \text{a.s..}$$
\end{thm}

\subsection{Our main results}

Motivated by Takei \cite{TMA}, we consider the LERRW $\boldsymbol{X}$  on $ \mathbb{Z}_{+}$ with  initial weights
\begin{align}\label{2.1}
w_0(x)=
\left\{
\begin{array}{cl}
1 & (x=0, 1),\\
x^{\alpha}\ln^{\beta}x & (x
=2,3,\cdots).
\end{array}
\right.
\end{align}
By Theorem \ref{thm-2.1}, $ \boldsymbol{X}$  is recurrent a.s. if and only if
$\alpha < 1$, or $\alpha=1$ and  $\beta \leq 1$.
Note that when $\beta=0$,
the initial weights  reduce to the initial weights of  \cite{TMA}. In the sequel, we
assume $\beta\neq 0$.
Our main results are Theorems \ref{thm-2.5} and \ref{thm-2.6} below.  Let
\begin{align}\label{2.2}
K(\alpha,\beta, \delta)=
\left\{
\begin{array}{cl}
\frac{1-\alpha}{2\delta} & (\alpha<0), \\
\frac{1}{2\delta}& (\alpha=0, \beta < 0), \\
\frac{1}{\delta} & (\alpha=0, \beta > 0),\\
\frac{1-\alpha}{\delta} & (0<\alpha<1).
\end{array}
\right.
\end{align}

\begin{thm}\label{thm-2.5}
  Assume that  $\alpha<1$ and $\delta>0$.
The LERRW  $\boldsymbol{X}$ on $\mathbb{Z}_+$  with initial weights \eqref{2.1} and  reinforcement parameter  $\delta$ satisfies that
\begin{align*}
\limsup_{n \to \infty} \frac{X_n}{(K(\alpha, \beta, \delta) (1-\alpha)^{-\beta}\ln n\,(\ln \ln n)^{\beta})^{1 / (1 - \alpha)}} =1\quad \text{a.s.,}
\end{align*}
where $K(\alpha,\beta, \delta)$ is defined as in \eqref{2.2}.
\end{thm}

\begin{thm}\label{thm-2.6}
Assume that $ \alpha=1$, $\beta \leq 1$ and  $\delta>0$, and consider the LERRW $\boldsymbol{X}$   on $\mathbb{Z}_+$ with initial weights \eqref{2.1} and reinforcement parameter  $\delta$.\\
\indent (i) If $\beta<0$,  then for any $\varepsilon > 0$,
$$ \limsup_{n\to \infty}\frac{X_{n}}{\exp \left\{ \left(\ln n \right)^{(1-\varepsilon)/(1-\beta)} \right\}}=+\infty
\quad \text{and} \quad  \lim_{n \to \infty} \frac{X_{n}}{ \exp \left\{  \left (\ln n \right)^{(1+ \varepsilon)/(1-\beta)}     \right \}  }=0 \quad  \text {a.s.}.$$
\indent (ii) If  $0<\beta \leq 1$,  then for any $\varepsilon >0 $,
$$\limsup _{n \rightarrow \infty} \frac{X_{n}}{n^{(1-\varepsilon) / 2}}=+\infty
\quad \text { and } \quad \lim _{n \rightarrow \infty} \frac{X_{n}}{n^{(1+\varepsilon) / 2}}=0 \quad \text { a.s.. }$$
\end{thm}

\begin{rmk}\label{rem-2.7}
 Comparing Theorem \ref{thm-2.2} and Theorem \ref{thm-2.5},  Theorem \ref{thm-2.3} and Theorem \ref{thm-2.6}, respectively, we see that,
any $\beta\neq 0$ changes the limit behavior of the LERRW in the case $\alpha<1$, and only $\beta<0$ changes the limit behavior of the LERRW in the case $\alpha=1$.
\end{rmk}

To see the effect of the reinforcement, we also prove the following  result for the un-reinforced case.

\begin{thm}\label{thm-2.9}
Consider the un-reinforced random walk $\boldsymbol{X}$ on $\mathbb{Z}_+$  with the initial weights \eqref{2.1}.\\
\indent (i) If $\alpha < -1$, there are two subcases:\\
\indent\quad (i.1) If $\beta>0$, then for any $\varepsilon\in (0,(1-\alpha)/\beta)$,
$$\limsup _{n \rightarrow \infty} \frac{X_{n}}{n^{1 /(1-\alpha)}}>0
\quad \text { and } \quad \lim _{n \rightarrow \infty} \frac{X_{n}}{\left\{n(\ln n)^{1+\varepsilon}\right\}^{1 /(1-\alpha-\varepsilon \beta )}}=0 \quad \text { a.s.. }$$
\indent\quad (i.2) If $\beta<0$, then for any $\varepsilon >0$,
$$\limsup _{n \rightarrow \infty} \frac{X_{n}}{n^{1 /(1-\alpha-\varepsilon \beta )}}=+\infty
\quad \text { and } \quad \lim _{n \rightarrow \infty} \frac{X_{n}}{\left\{n(\ln n)^{1+\varepsilon}\right\}^{1 /(1-\alpha)}}=0 \quad \text { a.s.. }$$
\indent (ii) If $\alpha = -1$, there are two subcases:\\
\indent\quad (ii.1) If $\beta < -1 $, then for any $\varepsilon > 0$,
$$\limsup _{n \rightarrow \infty} \frac{X_{n}}{n^{1/ (2-\varepsilon \beta )}}=+\infty
\quad \text { and } \quad \lim _{n \rightarrow \infty} \frac{X_{n}}{\left\{n(\ln n)^{1+\varepsilon}\right\}^{1 / 2}}=0 \quad \text { a.s.. }$$
\indent\quad\ (ii.2) If $\beta \geq -1$, then for any $\varepsilon > 0$,
$$\limsup _{n \rightarrow \infty} \frac{X_{n}}{n^{(1-\varepsilon) / 2}}=+\infty
\quad \text { and } \quad \lim _{n \rightarrow \infty} \frac{X_{n}}{\left\{n(\ln n)^{1+\varepsilon}\right\}^{1 / 2}}=0 \quad \text { a.s.. }$$
 \indent (iii)  If $-1 < \alpha <1$, or $\alpha=1$ and $\beta\leq 1$, then for any $\varepsilon > 0$,
$$\limsup _{n \rightarrow \infty} \frac{X_{n}}{n^{1 / 2}}>0
\quad \text { and } \quad \lim _{n \rightarrow \infty} \frac{X_{n}}{\left\{n(\ln n)^{1+\varepsilon}\right\}^{1 / 2}}=0 \quad \text { a.s.. }$$
\end{thm}

\section{Un-reinforced random walks on $\mathbb{Z}_+$}\setcounter{equation}{0}

 This section is heavily inspired by Sections 3 and 4 of Takei \cite{TMA}.
In Subsection \ref{ss:3.1}, we collect some results for  (un-reinforced) random walks on $\mathbb{Z}_+$ with general (initial) weights. Then in Subsection \ref{ss:3.2}, we apply the results of
Subsection \ref{ss:3.1} to prove Theorem \ref{thm-2.9}. The results of Subsection \ref{ss:3.1} will also be used in the proofs of Theorems \ref{thm-2.5} and \ref{thm-2.6}.

\subsection{Some useful results about random walks with general weights}\label{ss:3.1}

Let $(w_x: x\in \mathbb{Z}_+)$ be a family of positive numbers and let $w_{-1}=0$. Define
\begin{align}\label{3.1}
p_{x}=\frac{w_{x}}{w_{x-1}+w_{x}}, \quad \text { and } \quad q_{x}=\frac{w_{x-1}}{w_{x-1}+w_{x}} \quad \text { for } x \in \mathbb{Z}_{+},
\end{align}
and
\begin{align}\label{e:rsgamma}
\gamma_{0}:=1, \quad \text { and } \quad \gamma_{x}:=\prod_{i=1}^{x} \frac{q_{i}}{p_{i}} \quad \text { for } x \in \mathbb{N} .
\end{align}
 By (\ref{3.1}) and (\ref{e:rsgamma}),  we have $\gamma_x=\frac{w_0}{w_x}$ for all $x\in \mathbb{Z}_+$.

We will consider
a random walk
$\boldsymbol{Y}=(Y_n,n\geq 0)$ on $\mathbb{Z}_+$ with weights
 $(w_x: x\in \mathbb{Z}_+)$, i.e., a random walk which, when at $x\in \mathbb{Z}_+$, moves to $x+1$ with probability $p_x$ and to $x-1$ with probability $q_x$ in the next step. $ \gamma_{x}$ (resp. $w_x$)  can be interpreted as  the resistance (resp. conductance) of the edge $ \{x, x+1\}$. For any $i\in \mathbb{Z}_+$, we will use
$\mathbf{P}_i$ to denote the law of
$\boldsymbol{Y}$ when $Y_0=i$
and use $\mathbf{E}_i$ to denote the corresponding expectation.
Let
\begin{align}\label{e:rsharmonic}
h(0)=0, \quad h(x):=\sum_{i=0}^{x-1} \gamma_{i}, \quad x\in \mathbb{N}.
\end{align}
The effective resistance between the origin and infinity is $h(\infty)=\sum \limits_{i=0}^{\infty} \gamma_{i}$. Define $ \left\{\pi_{x}\right\}_{x \in \mathbb{Z}_{+}}$ by
\begin{align}\label{e:rspi}
\pi_{x}:=w_{x-1}+w_{x} \quad \text { for } x \in \mathbb{Z}_{+}.
\end{align}
 Note that
\begin{align}\label{3.2}
Z:=\sum_{i=0}^{\infty} \pi_{i}<+\infty \quad \text { if and only if } \quad \sum_{i=0}^{\infty} \frac{1}{\gamma_{i}}<+\infty.
\end{align}
\indent For $x\in\mathbb{Z}_+$, let
\begin{align}\label{e:rshittingtime}
\tau_{x}:=\inf \left\{n \geq 0: Y_{n}=x
\right\}
\end{align}
be the first entrance time of $x$.
Define for $x\in\mathbb{Z}_+$,
\begin{align}\label{3.3}
T(x):=\sum_{j=0}^{x-1} \pi_{j}\{h(x)-h(j)\}=\sum_{j=0}^{x-1} \pi_{j} \sum_{i=j}^{x-1} \gamma_{i}=\sum_{i=0}^{x-1} \gamma_{i} \sum_{j=0}^{i} \pi_{j}.
\end{align}
It follows from \cite[Proposition 2.20]{LP} that
$\mathbf{E}_{0}\left[\tau_{x}\right]=T(x).$

The following lemma follows immediately from \cite[Lemma 6.1.4, Theorem 2.8.1]{MPW}. (Although
\cite[Lemma 6.1.4]{MPW} is stated for random walks with random weights, its proof is for
fixed weights.)

\begin{lem}\label{lem-3.3}
 Let  $t_{1}(\cdot)$  be a nonnegative increasing function on $ \mathbb{Z}_{+}$ with $t_{1}(x) \rightarrow \infty$  as $ x \rightarrow \infty$. If
$T(x) \geq t_{1}(x)$  for all but finitely many $x \in \mathbb{Z}_{+},$
then for any $ \varepsilon>0$,   it holds   $\mathbf{P}_{0} \text {-a.s.}$ that
\begin{align*}
 Y_{n} \leq t_{1}^{-1}\left(2 n\{\ln (2 n)\}^{1+\varepsilon}\right) \quad \text{ for all but finitely many}  \ n.\end{align*}
\end{lem}

The following result is
 \cite[Lemma 4.3]{HMW}. Although
\cite[Lemma 4.3]{HMW} is stated for random walks with random weights, it is proved there for
fixed weights.

\begin{lem}\label{lem-3.4} Let  $t_{2}(\cdot)$  be a nonnegative increasing function on  $\mathbb{Z}_{+}$ with
$\sum_{x=1}^{\infty} \frac{t_{2}(x)}{t_{2}\left(x^{2}\right)}<\infty.$
If
$T(x) \leq t_{2}(x)$  for infinitely many $x \in \mathbb{Z}_{+}$,
then for any $ \varepsilon>0 $, it holds $\mathbf{P}_{0} \text {-a.s.}$  that
\begin{align*}
 Y_{n} \geq t_{2}^{-1}((1-\varepsilon) n) \quad  \text{for infinitely many }  n .
\end{align*}
\end{lem}

We now collect several useful bounds for $T(x)$.
\eqref{3.4} and \eqref{3.5} are contained in \cite[Lemma 3.5]{HMW}, \eqref{3.6} follows directly from
 \eqref{3.3}.

\begin{lem}\label{lem-3.5}
 For any  $x \in \mathbb{Z}_{+}$, we have
\begin{align}
&T(x) \geq h(x) \geq \max _{0 \leq i<x} \gamma_{i} \geq \gamma_{x-1},\label{3.4}\\
&T(x) \leq 2 x^{2}\left(\max _{0 \leq i<x} \gamma_{i}\right)\left(\max _{0 \leq j<x} \frac{1}{\gamma_{j}}\right).\label{3.5}
\end{align}
If  $Z=\sum \limits _{i=0}^{\infty} \pi_{i}<\infty$, then \eqref{3.5} can be improved as follows:
\begin{align}\label{3.6}
T(x) \leq Z h(x) \leq Z x\left(\max _{0 \leq i<x} \gamma_{i}\right).
\end{align}
\end{lem}

\subsection{Proof of Theorem \ref{thm-2.9}}\label{ss:3.2}

Now we return to the un-reinforced random walk $\boldsymbol{X}$ on $\mathbb{Z}_+$  with the initial weights \eqref{2.1} and apply the results of
Subsection \ref{ss:3.1}
to prove Theorem \ref{thm-2.9}.
Note that $\boldsymbol{X}$ is a particular case of $\boldsymbol{Y}$.
We will also use the following  infinite series version of l'H\^{o}pital's rule, due to Stolz and Ces\`{a}ro.

\begin{lem}\label{lem-3.6} {\rm (see e.g. Knopp \cite{K}, p. 34)}
 If a real sequence $\left\{a_{n}\right\}$ and a positive sequence $\left\{b_{n}\right\}$ satisfy
$\lim\limits_{n \rightarrow \infty} \frac{a_{n}}{b_{n}}=L \in \mathbb{R} \cup\{ \pm \infty\}$  and $ \sum_{n=1}^{\infty} b_{n}=+\infty$,
then
$\lim\limits_{n \rightarrow \infty} \frac{\sum_{k=1}^{n} a_{k}}{\sum_{k=1}^{n} b_{k}}=L.$
\end{lem}

\begin{lem}\label{lem-3.7}
\indent (i)
If $\rho>-1$ and  $L$ is a slowly varying sequence, then
\begin{align*}
\sum_{i=2}^{n} i^\rho L(i)\sim \frac{1}{\rho+1} n^{\rho+1}L(n).
\end{align*}

\indent (ii) For $\beta\leq 1$,
\begin{align*}
\sum_{i=2}^{n}\frac{1}{i(\ln i)^\beta}\sim
\begin{cases}
\dfrac{1}{1-\beta}(\ln n)^{1-\beta}, & \beta<1,\\
\ln\ln n, & \beta=1.
\end{cases}
\end{align*}
\end{lem}

\noindent{\it Proof.}
(i) This follows immediately from  the discrete Karamata's theorem \cite[Corollary 1.39]{GdH}.

\indent (ii) Let
\begin{align*}
f(x)=\frac{1}{x(\ln x)^\beta}, \quad x\geq 2.
\end{align*}
Then
\begin{align*}
\sum_{i=2}^{n}\frac{1}{i(\ln i)^\beta}\sim\int_2^n \frac{dx}{x(\ln x)^\beta}.
\end{align*}
If $\beta<1$, then
\begin{align*}
\int_2^n \frac{dx}{x(\ln x)^\beta}=\frac{1}{1-\beta}\left\{(\ln n)^{1-\beta}-(\ln 2)^{1-\beta}\right\}
\sim\frac{1}{1-\beta}(\ln n)^{1-\beta}.
\end{align*}
If $\beta=1$, then
\begin{align*}
\int_2^n \frac{dx}{x\ln x}=\ln\ln n-\ln\ln 2\sim\ln\ln n.
\end{align*}
The assertion now follows immediately.\hfill\fbox

\bigskip 

\noindent {\it  Proof of Theorem \ref{thm-2.9}.}
Recall that
$\beta\neq 0$.  Assume that $\alpha < 1$,  or
$\alpha=1$ and $\beta \leq 1$.
Notice that $ Z=\sum_{j=0}^{\infty} \pi_{j}<+\infty $ if and only if $\alpha<-1$, or
$\alpha= -1$ and $\beta < -1$.

 (i) For $\alpha < -1$, we deal with the two subcases $\beta>0$ and $\beta<0$ separately.\\
\indent (i.1) $\alpha < -1$ and $\beta>0$. In this case, we have
\begin{align*}
h(x)=\sum_{i=0}^{x-1}\gamma_{i}=2+\sum_{i=2}^{x-1}\frac{1}{i^{\alpha}\ln^{\beta}i}.
\end{align*}
It suffices to consider $\varepsilon\in (0,\min\{1,(1-\alpha)/\beta\})$. For such $\varepsilon$, there exists $N \geq 2$ such that $1 \leq \ln^{\beta}i \leq i^{\varepsilon \beta}$ for all $i \geq N$. Thus
\begin{align*}
\sum_{i=N}^{x-1}\frac{1}{i^{\alpha} \cdot i^{\varepsilon \beta}}\leq \sum_{i=N}^{x-1}\frac{1}{i^{\alpha}\ln^{\beta}i}\leq \sum_{i=N}^{x-1}\frac{1}{i^{\alpha}}.
\end{align*}
Since
\begin{align*}
\sum_{i=N}^{x-1}\frac{1}{i^{\alpha} \cdot i^{\varepsilon \beta}} \sim \frac{1}{1-\alpha-\varepsilon \beta}x^{1-\alpha-\varepsilon \beta}, \quad  \sum_{i=N}^{x-1}\frac{1}{i^{\alpha}} \sim \frac{1}{1-\alpha}x^{1-\alpha}\quad \mbox{as} \ x\to\infty,
\end{align*}
we have, for all sufficiently large $x$,
\begin{align*}
\frac{1-\varepsilon}{1-\alpha-\varepsilon \beta}x^{1-\alpha-\varepsilon \beta} \leq h(x)\leq \frac{1+\varepsilon}{1-\alpha}x^{1-\alpha}.
\end{align*}
Then by Lemma \ref{lem-3.5}, we get  that for such $\varepsilon$,
\begin{align*}
\frac{1-\varepsilon}{1-\alpha-\varepsilon \beta}x^{1-\alpha-\varepsilon \beta} \leq T(x)\leq \frac{(1+\varepsilon)Z}{1-\alpha}x^{1-\alpha} \quad \text{for all but finitely many } x.
\end{align*}
Hence, applying Lemma \ref{lem-3.3} with parameter $\varepsilon/2$, Lemma \ref{lem-3.4}, and using the fact that $ Z \geq \pi_{0}=1$, we get that
$\mathbf{P}_{0}$-a.s.,
\begin{align*}
X_n \leq \left(\frac{1-\alpha-\varepsilon\beta}{1-\varepsilon}\cdot 2n\{\ln(2n)\}^{1+\varepsilon/2}\right)^{1/(1-\alpha-\varepsilon\beta)} \quad \text{for all large } n
\end{align*}
and
\begin{align*}
X_n \geq \left\{\frac{1-\alpha}{(1+\varepsilon)Z}\cdot(1-\varepsilon)n\right\}^{1/(1-\alpha)} \quad \text{for infinitely many } n.
\end{align*}
The assertions of (i.1) now follow immediately.

\indent (i.2) $\alpha < -1$ and $\beta<0$. In this case, we have
\begin{align*}
h(x)=\sum_{i=0}^{x-1}\gamma_{i}=2+\sum_{i=2}^{x-1}\frac{1}{i^{\alpha}\ln^{\beta}i}.
\end{align*}
It suffices to consider $\varepsilon\in (0,1)$. For such $\varepsilon$,
there exists $N \geq 2$ such that $i^{\varepsilon \beta} \leq \ln^{\beta}i \leq 1$ for all $i\geq N$. Thus
\begin{align*}
\sum_{i=N}^{x-1}\frac{1}{i^{\alpha}}\leq \sum_{i=N}^{x-1}\frac{1}{i^{\alpha}\ln^{\beta}i}\leq \sum_{i=N}^{x-1}\frac{1}{i^{\alpha}\cdot i^{\varepsilon \beta}}.
\end{align*}
Since
\begin{align*}
\sum_{i=N}^{x-1}\frac{1}{i^{\alpha}} \sim \frac{1}{1-\alpha}x^{1-\alpha}, \quad  \sum_{i=N}^{x-1}\frac{1}{i^{\alpha}\cdot i^{\varepsilon \beta}} \sim \frac{1}{1-\alpha-\varepsilon \beta}x^{1-\alpha-\varepsilon \beta}\quad \mbox{as} \ x\to\infty,
\end{align*}
we have, for all sufficiently large $x$,
\begin{align*}
\frac{1-\varepsilon}{1-\alpha}x^{1-\alpha} \leq h(x)\leq \frac{1+\varepsilon}{1-\alpha-\varepsilon \beta}x^{1-\alpha-\varepsilon \beta}.
\end{align*}
Then by Lemma \ref{lem-3.5}, we get that
for $\varepsilon\in (0,1)$,
\begin{align*}
\frac{1-\varepsilon}{1-\alpha}x^{1-\alpha} \leq T(x)\leq \frac{(1+\varepsilon)Z}{1-\alpha-\varepsilon \beta}x^{1-\alpha-\varepsilon \beta} \quad \text{for all but finitely many } x.
\end{align*}
Hence, applying Lemma \ref{lem-3.3} with parameter $\varepsilon/2$, Lemma \ref{lem-3.4}, and using the fact that $ Z \geq \pi_{0}=1$, we get that
$\mathbf{P}_{0}$-a.s.,
\begin{align*}
X_n \leq \left(\frac{1-\alpha}{1-\varepsilon}\cdot 2n\{\ln(2n)\}^{1+\varepsilon/2}\right)^{1/(1-\alpha)} \quad \text{for all large } n
\end{align*}
and
\begin{align*}
X_n \geq \left\{\frac{1-\alpha-\varepsilon\beta}{(1+\varepsilon)Z}\cdot(1-\varepsilon)n\right\}^{1/(1-\alpha-\varepsilon\beta)} \quad \text{for infinitely many } n.
\end{align*}
By taking a smaller parameter if necessary, we obtain the conclusion of (i.2).

(ii) For $\alpha = -1$, we deal with the two subcases $\beta<-1$ and $\beta\geq -1(\beta\neq 0)$
separately. It suffices to consider $\varepsilon\in(0,1)$.\\
\indent (ii.1)
$\alpha = -1$ and $\beta < -1$. In this case, we have
\begin{align*}
h(x) = 2 + \sum_{i=2}^{x-1} \frac{1}{i^{-1} \ln^{\beta} i}.
\end{align*}
For $\varepsilon\in(0,1)$, there exists $N \geq 2$ such that $i^{\varepsilon \beta} \leq \ln^{\beta}i \leq 1$ for all $i \geq N$. Thus
\begin{align*}
\sum_{i=N}^{x-1}\frac{1}{i^{-1}}\leq \sum_{i=N}^{x-1}\frac{1}{i^{-1}\ln^{\beta}i}\leq \sum_{i=N}^{x-1}\frac{1}{i^{-1}\cdot i^{\varepsilon \beta}}.
\end{align*}
Since
\begin{align*}
\sum_{i=N}^{x-1}\frac{1}{i^{-1}}\sim \frac{1}{2}x^{2}, \quad
\sum_{i=N}^{x-1}\frac{1}{i^{-1}\cdot i^{\varepsilon \beta}}\sim \frac{1}{2-\varepsilon \beta}x^{2-\varepsilon \beta}\quad \mbox{as}\ x\to \infty,
\end{align*}
we have, for all sufficiently large $x$,
\begin{align*}
\frac{1-\varepsilon}{2}x^{2} \leq h(x) \leq \frac{1+\varepsilon}{2-\varepsilon \beta}x^{2-\varepsilon \beta}.
\end{align*}
Then by Lemma \ref{lem-3.5}, we get that for $\varepsilon\in(0,1)$,
\begin{align*}
\frac{1-\varepsilon}{2}x^{2} \leq T(x) \leq \frac{(1+\varepsilon)Z}{2-\varepsilon \beta}x^{2-\varepsilon \beta}  \quad  \text{for all but finitely many} \ x .
\end{align*}
The rest of the proof is similar to that of (i).\\
\indent  (ii.2)
$\alpha=-1$ and $\beta\geq -1(\beta\neq 0)$. We decompose this case into two subcases:
$\alpha=-1$ and $\beta=-1$, and $\alpha=-1$ and $\beta>-1$.\\
\indent  (ii.2.1) $\alpha = -1$ and $\beta=-1$. In this case, we have
\begin{align*}
\gamma_i  \sum_{j=0}^{i} \pi_j = i \ln i \left(4 + 2 \sum_{j=2}^{i} \frac{1}{j \ln j} - \frac{1}{i \ln i}\right)
\sim2 i \ln i \ln \ln i \quad \text{as } i \to \infty.
\end{align*}
By Lemma \ref{lem-3.6} and Lemma \ref{lem-3.7}(i),
\begin{align}\label{4.30}
T(x)=\sum_{i=0}^{x-1} \gamma_{i} \sum_{j=0}^{i} \pi_{j} \sim \sum_{i=2}^{x-1}2i\ln i \ln\ln i
\sim x^{2}\ln x \ln \ln x \quad \text{as} \quad x\to\infty.
\end{align}
Consequently, for any $\varepsilon\in(0,1)$,
\begin{align*}
  (1-\varepsilon) x^{2} \leq T(x) \leq(1+\varepsilon) x^{2+\varepsilon} \quad  \text{for all but finitely many} \  x .
\end{align*}
The rest of the proof is similar to that of (i).\\
\indent (ii.2.2)
$\alpha = -1$ and $\beta> -1$. In this case, we have
\begin{align*}
\gamma_i  \sum_{j=0}^{i} \pi_j = i (\ln i)^{-\beta}  \left(4 + 2 \sum_{j=2}^{i} \frac{1}{j (\ln j)^{-\beta} } - \frac{1}{i (\ln i)^{-\beta}}\right)
\sim\frac{2}{\beta+1} i \ln i \quad \text{as } i \to \infty.
\end{align*}
By Lemma \ref{lem-3.6} and Lemma \ref{lem-3.7}(i),
 \begin{align*}
T(x)=\sum_{i=0}^{x-1} \gamma_{i} \sum_{j=0}^{i} \pi_{j} \sim \frac{1}{\beta+1}x^{2} \ln x \quad \text { as } x \rightarrow \infty .
\end{align*}
Hence, for any $\varepsilon\in(0,1)$,
 \begin{align*}
 \frac{1-\varepsilon}{\beta+1} x^{2} \leq T(x) \leq \frac{1+\varepsilon}{\beta+1} x^{2+\varepsilon} \quad  \text{for all but finitely many} \  x .
 \end{align*}
The rest of the proof is similar to that of (i).

(iii) For $-1<\alpha<1$, or $\alpha=1$ and $\beta\le1$, we deal with the two subcases separately.\\
\indent (iii.1) $-1 < \alpha <1$.
In this case, by Lemma \ref{lem-3.7}(i), we have
\begin{align*}
\gamma_i  \sum_{j=0}^{i} \pi_j& = i^{-\alpha} (\ln i)^{-\beta}  \left(4 + 2 \sum_{j=2}^{i} \frac{1}{j^{-\alpha} (\ln j)^{-\beta}} - \frac{1}{i^{-\alpha} (\ln i)^{-\beta}}\right)\\
& \sim i^{-\alpha} (\ln i)^{-\beta}\cdot \frac{2}{\alpha+1}i^{1+\alpha}(\ln i)^{\beta}\\
&=\frac{2}{\alpha+1} i  \quad\quad \text{as } i \to \infty,
\end{align*}
and thus by Lemma \ref{lem-3.6},
\begin{align*}
T(x)=\sum_{i=0}^{x-1} \gamma_{i} \sum_{j=0}^{i} \pi_{j} \sim \frac{1}{\alpha+1}x^{2} \quad \text { as } x \rightarrow \infty.
\end{align*}
It follows that for any $\varepsilon\in(0,1)$,
\begin{align*}
 \frac{1-\varepsilon}{\alpha+1} x^{2} \leq T(x) \leq \frac{1+\varepsilon}{\alpha+1} x^{2} \quad  \text{for all but finitely many} \  x .
 \end{align*}
The rest of the proof is similar to that of (i).

\indent (iii.2) $\alpha=1$ and $\beta\leq1$.
In this case, by Lemma \ref{lem-3.7}(i), we have
\begin{align*}
\gamma_i  \sum_{j=0}^{i} \pi_j
&=i^{-1}(\ln i)^{-\beta}\left(4+2\sum_{j=2}^{i}\frac{1}{j^{-1}(\ln j)^{-\beta}}-\frac{1}{i^{-1}(\ln i)^{-\beta}}\right)\\
&\sim i^{-1}(\ln i)^{-\beta}\cdot i^2(\ln i)^{\beta}\\
&=i\quad \text{as } i\to\infty.
\end{align*}
By Lemma \ref{lem-3.6},
\begin{align*}
T(x)=\sum_{i=0}^{x-1}\gamma_i\sum_{j=0}^{i}\pi_j\sim \frac{1}{2}x^2\quad \text{as } x\to\infty.
\end{align*}
The rest of the proof is similar to that of (i).\hfill\fbox{}

\section{Proofs of Theorems \ref{thm-2.5} and \ref{thm-2.6}}\label{s:main}\setcounter{equation}{0}

In this section, we will give the proofs of Theorems \ref{thm-2.5} and \ref{thm-2.6}. We assume throughout this section that $\delta>0$ and  $(w_0(i): i\in \mathbb{Z}_+)$ are the initial weights given in \eqref{2.1}.
As in \cite{TMA}, we will prove Theorems \ref{thm-2.5} and \ref{thm-2.6} with the help of a random walk in random environment.
Let  $p_{0}=1$ and $\delta>0$.  Suppose that  $\left\{p_{i}(\omega)\right\}_{i \in \mathbb{N}}$  is a sequence of independent random variables
with $p_{i}$ being a $\mathrm{Beta}\left(\frac{w_0(i)}{2 \delta}, \frac{w_0(i-1)+\delta}{2 \delta}\right)$ random variable,
that is,
a random variable with density
\begin{align*}
B\left(\frac{w_0(i)}{2 \delta}, \frac{w_0(i-1)+\delta}{2 \delta}\right)^{-1}u^{\frac{w_0(i)}{2 \delta}-1}(1-u)^{\frac{w_0(i-1)+\delta}{2 \delta}-1}1_{(0, 1)}(u),
\end{align*}
where
\begin{align*}
B\left(\frac{w_0(i)}{2 \delta}, \frac{w_0(i-1)+\delta}{2 \delta}\right)=\int_{0}^{1} t^{\frac{w_0(i)}{2 \delta}-1}(1-t)^{\frac{w_0(i-1)+\delta}{2 \delta}-1} d t.
\end{align*}
We use $\mathbb{P}$ to denote the law of the random environment.
The expectation and variance under  $\mathbb{P}$  are denoted by  $\mathbb{E}[\cdot]$  and  $\mathbb{V}[\cdot]$, respectively.

Using the $ \left\{p_{i}(\omega)\right\}_{i \in \mathbb{Z}_{+}}$ above, we can define
a Markov chain
$\boldsymbol{Y}=\left\{Y_{n}\right\} $ on $\mathbb{Z}_{+}$
starting from $i_0\in \mathbb{Z}_+$ as follows:
$ \mathbf{P}_{i_{0}}^{\omega}\left(Y_{0}=i_{0}\right)=1$ and for any $n\ge 0$,
\begin{align*} \begin{cases}
\mathbf{P}_{i_{0}}^{\omega}\left(Y_{n+1}=i+1 \mid Y_{n}=i\right)=p_{i}(\omega),\\
\mathbf{P}_{i_{0}}^{\omega}\left(Y_{n+1}=i-1 \mid Y_{n}=i\right)=q_{i}(\omega):=1-p_{i}(\omega).
\end{cases}
\end{align*}
According to \cite[Section 3]{PR}, this RWRE is related to our LERRW as follows:

\begin{lem}\label{lem-3.1}
 For any  $n \geq 0$  and any $ i_{0}, i_{1}, \cdots, i_{n} \in \mathbb{Z}_{+} $, it holds that
\begin{align*}
P\left(X_{1}=i_{1}, \ldots, X_{n}=i_{n} \mid X_{0}=i_{0}\right)=\mathbb{E}\left[\mathbf{P}_{i_{0}}^{\omega}\left(Y_{1}=i_{1}, \ldots, Y_{n}=i_{n}\right)\right].
\end{align*}
\end{lem}
\indent  Define $\{\gamma_x(\omega)\}_{x\in \mathbb{Z}_+}$ by
 \begin{align*}
\gamma_{0}(\omega):=1, \quad \text { and } \quad \gamma_{x}(\omega):=\prod_{i=1}^{x} \frac{q_{i}(\omega)}{p_{i}(\omega)} \quad \text { for } x \in \mathbb{N} .
\end{align*}
Setting  $w(x, \omega):=1 / \gamma_{x}(\omega)$, we get
\begin{align}\label{3.1-1}
p_{x}(\omega)=\frac{w(x, \omega)}{w(x-1, \omega)+w(x, \omega)}, \quad q_{x}(\omega)=\frac{w(x-1, \omega)}{w(x-1, \omega)+w(x, \omega)} \quad \text { for } x \in \mathbb{Z}_{+},
\end{align}
where  $w(-1, \omega):=0$. Thus for any $\omega$, under $\mathbf{P}^\omega_{i_{0}}$,
$\boldsymbol{Y}$ is a random walk on $\mathbb{Z}_+$ with weights
$(w(x, \omega): x\in \mathbb{Z}_+)$. Thus, for any fixed $\omega$, we can apply the results of
Subsection \ref{ss:3.1}.

For simplicity, we will sometimes suppress the superscript $\omega$ in the remainder of this section.
We need a result
on the asymptotic behavior of $\ln \gamma_x$ as $x\to \infty$.
Let
\begin{align*}
S_{x}:=\ln \gamma_{x}=\ln\prod_{i=1}^{x} \frac{q_{i}}{p_{i}}=\sum_{i=1}^{x} \ln \frac{1-p_{i}}{p_{i}} \quad \text { for } x \in \mathbb{N} .
\end{align*}

\begin{pro}\label{pro-5.1}
\indent (i) If  $\alpha<1$, then
\begin{align*}
\lim _{x \rightarrow \infty} \frac{S_{x}}{x^{1-\alpha}(\ln x)^{-\beta}}=\frac{1}{K(\alpha, \beta, \delta)} \quad \mathbb{P} \text {-a.e. } \omega,
\end{align*}
where $K(\alpha, \beta, \delta)$  is defined by \eqref{2.2}.\\
\indent (ii) If  $\alpha=1, \beta<0$, then
\begin{align*}
\lim _{x \rightarrow \infty} \frac{S_{x}}{(\ln x)^{1-\beta}}=\frac{\delta}{1-\beta} \quad \mathbb{P} \text {-a.e. } \omega.
\end{align*}
\indent (iii) If $\alpha=1, 0<\beta \leq 1$, then
\begin{align*}
\lim _{x \rightarrow \infty} \frac{S_{x}}{ \ln x}=-1\quad \mathbb{P} \text {-a.e. } \omega.
\end{align*}
\end{pro}

The proof of  Proposition \ref{pro-5.1} is
pretty long and we postpone it
to the next section.
We now use it to prove our main results.
By Proposition \ref{pro-5.1}, we have
\begin{align}\label{4.1}
Z^{\omega}=\sum_{x=0}^{\infty} \pi_{x}(\omega)<+\infty  \quad  \text{if} \ &\alpha<1, \ \text{or}\  \alpha=1\ \text{and}\ \beta < 0,
\end{align}
where $\pi_x(\omega)=w(x-1, \omega)+w(x, \omega)$.

\smallskip

\noindent {\it Proof of Theorem \ref{thm-2.5}.}
 We fix $\alpha<1$.  We shall use the following elementary inversion estimate. For $c>0$, let
\begin{equation*}
t_c(x):=\exp\left\{c x^{1-\alpha}(\ln x)^{-\beta}\right\}.
\end{equation*}
Then $t_c$ is eventually increasing. If $u=t_c(x)$, then
\begin{align*}
\ln u&=c x^{1-\alpha}(\ln x)^{-\beta},\\
\ln\ln u&\sim (1-\alpha)\ln x.
\end{align*}
Hence
\begin{equation*}
x^{1-\alpha}\sim c^{-1}(1-\alpha)^{-\beta}\ln u\, (\ln\ln u)^{\beta}.
\end{equation*}
Therefore, as $u\to\infty$,
\begin{equation}\label{eq:inv}
t_c^{-1}(u)\sim \left\{c^{-1}(1-\alpha)^{-\beta}\ln u\, (\ln\ln u)^{\beta}\right\}^{1/(1-\alpha)}.
\end{equation}
Moreover, $t_c$ satisfies the summability condition in Lemma \ref{lem-3.4}.
By \eqref{3.4} and Proposition \ref{pro-5.1} (i),
for any  $\varepsilon \in (0,1)$,
\begin{equation*}
T^{\omega}(x) \geq \gamma_{x-1} \geq \exp \left(\frac{1-\varepsilon}{K(\alpha, \beta, \delta)} x^{1-\alpha}(\ln x)^{-\beta}\right) \quad \text { for all large } x.
\end{equation*}
Then by Lemma \ref{lem-3.3}, $\mathbb{P}$-a.e.  $\omega$  and $ \mathbf{P}_{0}^{\omega}$-a.s.,
using \eqref{eq:inv} with $c=(1-\varepsilon)/K(\alpha,\beta,\delta)$ yields
\begin{equation*}
Y_n\leq t_c^{-1}\left(2n\{\ln(2n)\}^{1+\varepsilon}\right)\sim \left\{\frac{K(\alpha,\beta,\delta)}{1-\varepsilon}(1-\alpha)^{-\beta}\ln n\, (\ln\ln n)^{\beta}\right\}^{1/(1-\alpha)} \quad \text { for all large } n,
\end{equation*}
which implies
\begin{equation*}
\limsup_{n \rightarrow \infty} \frac{Y_{n}}{\left\{K(\alpha,\beta,\delta)(1-\alpha )^{-\beta} \ln n\,(\ln \ln n)^{\beta}\right\}^{1 /(1-\alpha)}} \leq (1-\varepsilon)^{-1/(1-\alpha)}.
\end{equation*}
By Lemma \ref{lem-3.1}, letting $\varepsilon\downarrow0$, we have
\begin{equation*}
\limsup_{n \rightarrow \infty} \frac{X_{n}}{\left\{K(\alpha,\beta,\delta)(1-\alpha )^{-\beta} \ln n\,(\ln \ln n)^{\beta}\right\}^{1 /(1-\alpha)}} \leq 1\quad P \text {-a.s.}.
\end{equation*}
\indent Next, we consider the lower bound. Fix an arbitrary $\varepsilon > 0$. By Proposition \ref{pro-5.1}
(i), we get
\begin{equation*}
\max_{0 \leq i \leq x} \gamma_i \leq \exp\left( \frac{1 + \varepsilon/2}{K(\alpha, \beta, \delta)} x^{1 - \alpha} (\ln x)^{-\beta} \right) \quad \text{for all large } x.
\end{equation*}
By \eqref{4.1},
\begin{equation*}
Z^{\omega}x \leq \exp\left( \frac{\varepsilon/2}{K(\alpha, \beta, \delta)} x^{1 - \alpha} (\ln x)^{-\beta} \right) \quad \text{for all large } x.
\end{equation*}
It follows from \eqref{3.6} that
\begin{equation*}
T^{\omega }(x) \leq Z^{\omega}x\left(\max_{0 \leq i \leq x} \gamma_i \right) \leq \exp\left( \frac{1 + \varepsilon}{K(\alpha, \beta, \delta)} x^{1 - \alpha} (\ln x)^{-\beta} \right) \quad \text{for all large } x.
\end{equation*}
Then by Lemma \ref{lem-3.4}, $\mathbb{P}$-a.e.  $\omega$  and $ \mathbf{P}_{0}^{\omega}$-a.s.,
using \eqref{eq:inv} with $c=(1+\varepsilon)/K(\alpha,\beta,\delta)$ yields
\begin{equation*}
Y_n\geq t_c^{-1}((1-\varepsilon)n)\sim \left\{\frac{K(\alpha,\beta,\delta)}{1+\varepsilon}(1-\alpha)^{-\beta}\ln n\,(\ln\ln n)^{\beta}\right\}^{1/(1-\alpha)} \quad \text{for infinitely many} \ n,
\end{equation*}
which implies
\begin{equation*}
\limsup_{n \to \infty} \frac{Y_n}{\left\{K(\alpha, \beta, \delta) (1-\alpha)^{-\beta}\ln n\,(\ln \ln n)^{\beta}\right\}^{1 / (1 - \alpha)}} \geq (1+\varepsilon )^{-1/(1-\alpha)}.
\end{equation*}
By Lemma \ref{lem-3.1}, letting $\varepsilon\downarrow0$, we have
\begin{equation*}
\limsup_{n \to \infty} \frac{X_n}{\left\{K(\alpha, \beta, \delta) (1-\alpha)^{-\beta}\ln n\,(\ln \ln n)^{\beta}\right\}^{1 / (1 - \alpha)}} \geq 1\quad P \text {-a.s.}.
\end{equation*}
The proof is complete.\hfill\fbox

\smallskip

\noindent{\it Proof of Theorem \ref{thm-2.6}.}
We deal with the two subcases $\beta<0$ and $0<\beta\leq1$ separately.\\
\indent  (i)
$\alpha=1$ and $\beta<0$.
By \eqref{3.4} and Proposition \ref{pro-5.1}(ii),
for any $\varepsilon \in (0, \delta)$,
\begin{align*}
T^{\omega}(x) \geq \gamma_{x-1} \geq \exp\left(\frac{\delta-\varepsilon}{1-\beta}(\ln x)^{1-\beta}\right)\quad \quad \text{for all large } x.
\end{align*}
By Lemma \ref{lem-3.3}, $\mathbb{P}$-a.e.  $\omega$  and  $\mathbf{P}_{0}^{\omega}$-a.s.,
\begin{align*}
Y_{n} \leq  \exp   \left\{  \left [ \frac{1-\beta}{\delta-\varepsilon}\ln\left ( 2n( \ln (2 n))^{1+\varepsilon}\right ) \right ]^{1/(1-\beta)} \right \}   \quad \quad \text{for all large } n,
\end{align*}
which together with Lemma \ref{lem-3.1} implies
\begin{align*}
\lim_{n \to \infty} \frac{X_{n}}{ \exp \left\{  \left (\ln n \right)^{(1+ \varepsilon)/(1-\beta)} \right \}  }=0 \quad P \text {-a.s.}.
\end{align*}
\indent Now we turn to the lower bound.
For any $\varepsilon>0$,
by Proposition \ref{pro-5.1}(ii), we get that
\begin{align*}
\max_{0 \leq i < x} \gamma_{i} \leq \exp\left(\frac{\delta + \varepsilon / 2}{1-\beta}(\ln x)^{1-\beta}\right)\quad \quad \text{for all large } x.
\end{align*}
\indent By \eqref{4.1},
\begin{align*}
Z^{\omega}x \leq \exp\left(   \frac{\varepsilon / 2}{1-\beta} (\ln x)^{1-\beta}  \right)\quad \quad \text{for all large } x.
\end{align*}
It follows from \eqref{3.6} that
\begin{align*}
T^{\omega}(x) \leq \exp \left(   \frac{\delta + \varepsilon }{1-\beta} (\ln x)^{1-\beta} \right) \quad \quad \text{for all large } x.
\end{align*}
By Lemma \ref{lem-3.4}, $\mathbb{P}$-a.e.  $\omega$  and  $\mathbf{P}_{0}^{\omega}$-a.s.,
\begin{align*}
Y_{n} \geq \exp \left\{ \left [ \frac{1-\beta}{\delta+\varepsilon} \ln ((1-\varepsilon)n)  \right ]^{1/(1-\beta)} \right \}  \quad  \text{for infinitely many} \ n,
\end{align*}
which together with  Lemma \ref{lem-3.1} implies that
\begin{align*}
\limsup_{n\to \infty}\frac{X_{n}}{\exp \left\{ \left(\ln n \right)^{(1-\varepsilon)/(1-\beta)} \right\}}=+\infty \quad P \text {-a.s.. }
\end{align*}
\indent (ii)
$\alpha=1$ and $0<\beta\leq1$.
Proposition \ref{pro-5.1}(iii) implies that for any $\varepsilon > 0$,
\begin{align*}
\gamma_{i} \geq \exp( (-1- \varepsilon) \ln i ) = i^{-1-\varepsilon} \quad \text{for all large} \ i
\end{align*}
and
\begin{align*}
\gamma_{i} \leq \exp( (-1+ \varepsilon) \ln i ) = i^{-1+ \varepsilon} \quad \text{for all large} \ i.
\end{align*}
Notice that
\begin{align*}
T^{\omega}(x)=\sum_{i=0}^{x-1}\gamma_{i}\sum_{j=0}^{i}\pi_{j}=1+\sum_{i=1}^{x-1}\gamma_{i}\left(2+2\sum_{j=1}^{i}\frac{1}{\gamma_{j}}-\frac{1}{\gamma_{i}}\right)
\end{align*}
and
\begin{align*}
\sum_{i=1}^{x-1}i^{-1\pm \varepsilon}\sum_{j=1}^{i}j^{1 \pm \varepsilon} \sim \frac{1}{(2 \pm \varepsilon)(2\pm 2\varepsilon)}x^{2\pm2\varepsilon} \ \text{as} \ x \rightarrow \infty.
\end{align*}
\indent For any $\varepsilon > 0$, we have
\begin{align*}
T^{\omega}(x) \leq \frac{1+\varepsilon }{(2+\varepsilon)(2+2\varepsilon)}x^{2+2\varepsilon} \quad \text{for all large} \ x.
\end{align*}
\indent  On the other hand, for any $\varepsilon \in(0,1)$,
\begin{align*}
 T^{\omega}(x) \geq \frac{1-\varepsilon}{(2-\varepsilon)(2-2 \varepsilon)} x^{2-2 \varepsilon} \quad \text{ for all large}  \ x .
\end{align*}
The rest of the proof is similar to that of Theorem \ref{thm-2.9} (ii.2).
\hfill\fbox

\section{ Proof of Proposition \ref{pro-5.1} }\setcounter{equation}{0}

To derive the results of Proposition \ref{pro-5.1}, we need several lemmas. By (4.13) and (4.15) of \cite{TM}, for  $x \in \mathbb{N}$, we have
\begin{align}
\mathbb{E}\left[S_{x}\right] & =\sum_{i=1}^{x}\left\{\Psi\left(\frac{w_{0}(i-1)+\delta}{2 \delta}\right)-\Psi\left(\frac{w_{0}(i)}{2 \delta}\right)\right\}  \label{6.1} \\
& =\Psi\left(\frac{w_{0}(0)}{2 \delta}\right)-\Psi\left(\frac{w_{0}(x)}{2 \delta}\right)+\sum_{i=0}^{x-1}\left\{\Psi\left(\frac{w_{0}(i)}{2 \delta}+\frac{1}{2}\right)-\Psi\left(\frac{w_{0}(i)}{2 \delta}\right)\right\},\label{6.2}\\
\mathbb{V}\left[S_{x}\right] & =\sum_{i=1}^{x}\left\{\Psi^{\prime}\left(\frac{w_{0}(i-1)+\delta}{2 \delta}\right)+\Psi^{\prime}\left(\frac{w_{0}(i)}{2 \delta}\right)\right\},\label{6.3}
\end{align}
where $\Psi(z)=\Gamma'(z)/\Gamma(z)$ is the digamma function.
The following result is well-known.

\begin{lem}\label{lem-5.1} {\rm (\cite[Lemma 6.1]{TMA})}
\begin{align}\label{6.4}
\Psi^{\prime}(z) \sim \begin{cases}
\frac{1}{z} & (z \rightarrow \infty), \\
\frac{1}{z^{2}} & (z \rightarrow 0),
\end{cases}
\end{align}
and
\begin{align}\label{6.5}
\Psi\left(z+\frac{1}{2}\right)-\Psi(z) \sim \begin{cases}
\frac{1}{2 z} & (z \rightarrow \infty), \\
\frac{1}{z} & (z \rightarrow 0).
\end{cases}
\end{align}
\end{lem}
\indent   It follows from  \cite[(4.11)]{TM} that for any $y, z>0$,
\begin{align}\label{6.6}
\ln y-\ln z-\frac{1}{y} \leq \Psi(y)-\Psi(z) \leq \ln y-\ln z+\frac{1}{z}.
\end{align}

\begin{lem}\label{lem-5.2}
 Assume that $\delta>0$, $\beta \ne 0$, and that either $\alpha <1$, or $\alpha =1$ and $\beta \leq 1$. As $x\to\infty$,
\begin{align*}
\mathbb{E}[S_{x}] \sim \begin{cases}
 \frac{2\delta}{1-\alpha}x^{1-\alpha}(\ln x)^{-\beta}&(\alpha<0),\\
\delta x (\ln x)^{-\beta} & (\alpha=0, \beta > 0), \\
2\delta x (\ln x)^{-\beta} &  (\alpha=0, \beta < 0),  \\
\frac{\delta}{1-\alpha} x^{1-\alpha} (\ln x)^{-\beta}& (0<\alpha<1),\\
 \frac{\delta}{1-\beta}(\ln x)^{1-\beta}&(\alpha=1, \beta<0),\\
 -\ln x&(\alpha=1, 0<\beta\leq 1).
\end{cases}
\end{align*}
\end{lem}
\noindent{\it Proof.} By \eqref{6.6} and
\begin{align*}
\ln \frac{w_{0}(0)}{2 \delta}-\ln \frac{w_{0}(x)}{2 \delta}=\ln w_{0}(0)-\ln w_{0}(x)=-\ln [x^{\alpha}\ln^{\beta}x]=-\alpha \ln x-\beta \ln \ln x,
\end{align*}
the first term in \eqref{6.2} satisfies
\begin{align}\label{6.7}
-\alpha \ln x-\beta \ln \ln x-2 \delta \leq \Psi\left(\frac{w_{0}(0)}{2 \delta}\right)-\Psi\left(\frac{w_{0}(x)}{2 \delta}\right) \leq -\alpha \ln x-\beta \ln \ln x+\frac{2 \delta}{x^{\alpha}\ln^{\beta}x}.
\end{align}
\indent First, we assume that $0<\alpha \leq1$. By \eqref{6.7}, we have
\begin{align}\label{6.8}
\Psi\left(\frac{w_{0}(0)}{2 \delta}\right)-\Psi\left(\frac{w_{0}(x)}{2 \delta}\right) \sim -\alpha \ln x-\beta \ln \ln x \quad \text { as } x \rightarrow \infty.
\end{align}
Now  $w_{0}(x) \rightarrow \infty$  as  $x \rightarrow \infty$. By \eqref{6.5}, we have
\begin{align*}
\Psi\left(\frac{w_{0}(i)}{2 \delta}+\frac{1}{2}\right)-\Psi\left(\frac{w_{0}(i)}{2 \delta}\right) \sim \frac{\delta}{i^{\alpha}\ln^{\beta}i} \quad \text { as } i \rightarrow \infty,
\end{align*}
 which together with Lemma \ref{lem-3.7}(i) for $0<\alpha<1$ and Lemma \ref{lem-3.7}(ii) for $\alpha=1$ implies
\begin{align}\label{6.9}
&\sum_{i=0}^{x-1}\left\{\Psi\left(\frac{w_{0}(i)}{2 \delta}+\frac{1}{2}\right)-\Psi\left(\frac{w_{0}(i)}{2 \delta}\right)\right\} \nonumber \\
&\sim \delta \sum_{i=2}^{x-1}  \frac{1}{i^{\alpha}\ln^{\beta}i} \sim \begin{cases}
\frac{\delta}{1-\alpha} x^{1-\alpha} (\ln x)^{-\beta}&(0<\alpha<1),\\
\frac{\delta}{1-\beta}(\ln x)^{1-\beta}& (\alpha=1, \beta<1),\\
\delta \ln \ln x&(\alpha=1, \beta=1)
\end{cases}\quad \text { as } x \rightarrow \infty.
\end{align}
By \eqref{6.2}, \eqref{6.8} and \eqref{6.9}, we get the conclusions for $ 0<\alpha \leq 1$.\\
\indent  Next, we assume that $\alpha=0$. By \eqref{6.7}, we have
 \begin{align}\label{6.10}
 -\beta \ln \ln x-2 \delta \leq \Psi\left(\frac{w_{0}(0)}{2 \delta}\right)-\Psi\left(\frac{w_{0}(x)}{2 \delta}\right) \leq -\beta \ln \ln x+2 \delta(\ln x)^{-\beta}.
 \end{align}
 When $\beta >0$, $w_0(x)\to \infty$ as $x\to \infty$. When $\beta < 0$, $w_0(x) \to 0$ as $x \to \infty$. Then by \eqref{6.5},
 \begin{align*}
  \Psi\left(\frac{w_{0}(i)}{2 \delta}+\frac{1}{2}\right)-\Psi\left(\frac{w_{0}(i)}{2 \delta}\right) \sim \begin{cases}
  \frac{\delta}{\ln^{\beta}i}&(\beta > 0),\\
  \frac{2\delta}{\ln^{\beta}i}&(\beta<0)
  \end{cases} \quad \text { as } i \rightarrow \infty,
 \end{align*}
which together with Lemma \ref{lem-3.7}(i) implies that
 \begin{align*}
&\sum_{i=0}^{x-1}\left\{\Psi\left(\frac{w_{0}(i)}{2 \delta}+\frac{1}{2}\right)-\Psi\left(\frac{w_{0}(i)}{2 \delta}\right)\right\} \nonumber \\
&\sim \begin{cases}
\delta \sum_{i=2}^{x-1}\frac{1}{\ln^{\beta}i}\sim \delta x(\ln x)^{-\beta}& (\beta > 0), \\
2 \delta \sum_{i=2}^{x-1}\frac{1}{\ln^{\beta}i}\sim 2 \delta x(\ln x)^{-\beta}& (\beta < 0)
\end{cases}\quad \text { as } x \rightarrow \infty.
\end{align*}
Then by virtue of \eqref{6.2} and \eqref{6.10}, we obtain the conclusion for $\alpha=0$.\\
\indent  Finally, we assume that  $\alpha<0$. Now $w_{0}(x)\to 0$   as $x \rightarrow \infty$. By \eqref{6.5}, we have
\begin{align*}
\Psi\left(\frac{w_{0}(i)}{2 \delta}+\frac{1}{2}\right)-\Psi\left(\frac{w_{0}(i)}{2 \delta}\right) \sim \frac{2 \delta}{i^{\alpha}\ln^{\beta}i} \quad \text { as } i \rightarrow \infty,
\end{align*}
which together with Lemma \ref{lem-3.7}(i) implies that
\begin{align*}
\sum_{i=0}^{x-1}\left\{\Psi\left(\frac{w_{0}(i)}{2 \delta}+\frac{1}{2}\right)-\Psi\left(\frac{w_{0}(i)}{2 \delta}\right)\right\}\sim2\delta\sum_{i=2}^{x-1}\frac{1}{i^{\alpha}\ln^{\beta}i} \sim\frac{2\delta}{1-\alpha}x^{1-\alpha}(\ln x)^{-\beta}\quad\text{as}\ x\to\infty.
\end{align*}
 Then by virtue of \eqref{6.2} and  \eqref{6.7}, we obtain the conclusion for $\alpha<0$.\hfill\fbox

\smallskip

\begin{lem}\label{lem-5.3}
 Assume that $\delta>0$, $\beta \ne 0$, and that either $\alpha <1$, or $\alpha =1$ and $\beta \leq 1$. As $x\to\infty$,
\begin{align*}
\mathbb{V}[S_{x}] \sim \begin{cases}
 \frac{4\delta^2}{1-2\alpha}x^{1-2\alpha}(\ln x)^{-2\beta}&(\alpha<0),\\
 4 \delta x(\ln x)^{-\beta}& (\alpha=0, \beta > 0),\\
 4\delta^{2}x(\ln x)^{-2\beta}&(\alpha=0, \beta<0),\\
\frac{4\delta}{1-\alpha} x^{1-\alpha} (\ln x)^{-\beta}& (0<\alpha<1),\\
\frac{4\delta}{1-\beta}(\ln x)^{1-\beta}& (\alpha=1, \beta<1),\\
4\delta \ln \ln x&(\alpha=1, \beta=1).
\end{cases}
\end{align*}
\end{lem}
\noindent{\it Proof.} First we assume that $0<\alpha \leq 1$. Now  $w_{0}(x) \rightarrow \infty$  as $ x \rightarrow \infty$. By \eqref{6.4}, we get
\begin{align*}
\Psi^{\prime}\left(\frac{w_{0}(i)}{2 \delta}\right) & \sim \frac{2 \delta}{i^{\alpha}\ln^{\beta}i},
\end{align*}
and
\begin{align*}
\Psi^{\prime}\left(\frac{w_{0}(i-1)+\delta}{2 \delta}\right) & \sim \frac{2 \delta}{(i-1)^{\alpha}\ln^{\beta}(i-1)+\delta} \sim \frac{2 \delta}{i^{\alpha}\ln^{\beta}i} \quad \text { as } i \rightarrow \infty.
\end{align*}
Then by \eqref{6.3}  and Lemma \ref{lem-3.7}(i) for $0<\alpha<1$ and Lemma \ref{lem-3.7}(ii) for $\alpha=1$, we have
\begin{align*}
\mathbb{V}\left[S_{x}\right] &=\sum_{i=1}^{x}\left\{\Psi^{\prime}\left(\frac{w_{0}(i-1)+\delta}{2 \delta}\right)+\Psi^{\prime}\left(\frac{w_{0}(i)}{2 \delta}\right)\right\}\\
&\sim 4 \delta \sum_{i=2}^{x} \frac{1}{i^{\alpha}\ln^{\beta}i} \sim\begin{cases}
\frac{4\delta}{1-\alpha} x^{1-\alpha} (\ln x)^{-\beta}&(0<\alpha<1),\\
\frac{4\delta}{1-\beta}(\ln x)^{1-\beta}& (\alpha=1, \beta<1),\\
4\delta \ln \ln x&(\alpha=1, \beta=1)
\end{cases} \quad \text { as } x \rightarrow \infty.
\end{align*}
\indent Next, we assume that $\alpha=0$. When $\beta > 0$, $w_0(x) \to \infty$ as $x \to \infty$. When $\beta < 0$,
 $w_0(x) \to 0$ as $x \to \infty$. By \eqref{6.4}, we have
\begin{align*}
\Psi^{\prime}\left(\frac{w_{0}(i)}{2 \delta}\right) & \sim \begin{cases}
\frac{2\delta}{\ln^{\beta}i} & (\beta > 0),\\
\frac{4\delta^{2}}{\ln^{2\beta}i}&(\beta<0),
\end{cases}
\end{align*}
and
\begin{align*}
\Psi^{\prime}\left(\frac{w_{0}(i-1)+\delta}{2 \delta}\right)  \sim  \begin{cases}
\frac{2\delta}{\ln^{\beta}i} & (\beta > 0),\\
 \Psi^{\prime}(\frac{1}{2})>0 &(\beta<0)
\end{cases}
\quad \text { as } i \rightarrow \infty.
\end{align*}
Then by \eqref{6.3}  and Lemma \ref{lem-3.7}(i), we get
\begin{align*}
\mathbb{V}\left[S_{x}\right] &=\sum_{i=1}^{x}\left\{\Psi^{\prime}\left(\frac{w_{0}(i-1)+\delta}{2 \delta}\right)+\Psi^{\prime}\left(\frac{w_{0}(i)}{2 \delta}\right)\right\}\\
&\sim \begin{cases}
4\delta\sum_{i=2}^{x}\frac{1}{\ln^{\beta}i} \sim 4 \delta x(\ln x)^{-\beta} & (\beta > 0), \\
4\delta^{2}\sum_{i=2}^{x}\frac{1}{\ln^{2\beta}i} \sim 4 \delta^{2} x(\ln x)^{-2\beta} & (\beta < 0 )
\end{cases} \quad \text { as } x \rightarrow \infty.
\end{align*}
\indent Finally, we assume that  $\alpha<0$. Now $w_{0}(x)\to 0$  as $ x \rightarrow \infty$. By \eqref{6.4}, we have
\begin{align*}
\Psi^{\prime}\left(\frac{w_{0}(i)}{2 \delta}\right) & \sim\left(\frac{2 \delta}{i^{\alpha}\ln^{\beta}i}\right)^{2}=\frac{4 \delta^{2}}{i^{2 \alpha}\ln^{2\beta}i},
\end{align*}
and
\begin{align*}
\Psi^{\prime}\left(\frac{w_{0}(i-1)+\delta}{2 \delta}\right) & \rightarrow \Psi^{\prime}\left(\frac{1}{2}\right)>0 \quad \text { as } i \rightarrow \infty.
\end{align*}
Then by \eqref{6.3}  and Lemma \ref{lem-3.7}(i) with $\rho=-2\alpha$,
we get
\begin{align*}
\mathbb{V}\left[S_{x}\right] &=\sum_{i=1}^{x}\left\{\Psi^{\prime}\left(\frac{w_{0}(i-1)+\delta}{2 \delta}\right)+\Psi^{\prime}\left(\frac{w_{0}(i)}{2 \delta}\right)\right\}\\
&\sim \Psi^{\prime}\left(\frac{1}{2}\right) x+4 \delta^{2}\sum_{i=2}^{x} \frac{1}{i^{2 \alpha}\ln^{2\beta}i}\\
&\sim 4 \delta^{2}\sum_{i=2}^{x} \frac{1}{i^{2 \alpha}\ln^{2\beta}i} \sim\frac{4\delta^2}{1-2\alpha}x^{1-2\alpha}(\ln x)^{-2\beta}
 \quad \text { as } x \rightarrow \infty.
\end{align*}
The proof is complete.\hfill\fbox

\smallskip

The following is a consequence of Kolmogorov's strong law of large numbers  (see \cite{IK}, \cite{TMA}).

\begin{lem}\label{lem-5.4} Let  $\left\{\eta_{i}\right\}$ be a sequence of independent, square-integrable random variables, and  $S_{x}:=\sum_{i=1}^{x} \eta_{i}$. If  $\lim\limits_{x\to\infty}\mathbb{V}\left[S_{x}\right]=\infty$, then for any  $\epsilon >0$,
$\lim\limits_{x \rightarrow \infty} \frac{S_{x}-\mathbb{E}\left[S_{x}\right]}{\left(\mathbb{V}\left[S_{x}\right]\right)^{1 / 2+\epsilon }}=0$ a.s.. In particular, if
$\lim\limits_{x \rightarrow \infty} \frac{\left(\mathbb{V}\left[S_{x}\right]\right)^{1 / 2+\epsilon }}{\mathbb{E}\left[S_{x}\right]}=0$ for some $\epsilon >0,$
then
 $\lim\limits_{x \rightarrow \infty} \frac{S_{x}}{\mathbb{E}\left[S_{x}\right]}=1$  a.s..
\end{lem}

\smallskip

\noindent {\it  Proof of Proposition \ref{pro-5.1}}:  Combining Lemmas \ref{lem-5.2}, \ref{lem-5.3} and \ref{lem-5.4}, we immediately get the desired assertion.  \qed

\bigskip

\noindent {\bf\large Acknowledgments}\
 This work was supported by the National Natural Science Foundation of China (Grant Nos. 12171335, 12471139) and the Simons Foundation (\#960480).

\end{document}